\documentclass[twoside,a4paper,12pt,centertags]{amsart}
\usepackage[T1]{fontenc}
\usepackage[utf8]{inputenc}
\usepackage{subfiles}
\makeatletter
\usepackage[section]{placeins}
\AtBeginDocument{%
	\expandafter\renewcommand\expandafter\subsection\expandafter{%
		\expandafter\@fb@secFB\subsection
	}%
}
\makeatother
\makeatletter
\@namedef{subjclassname@2020}{%
	\textup{2020} Mathematics Subject Classification}
\makeatother
\usepackage{amsfonts,amssymb,amsmath,amsthm}
\usepackage{mathrsfs}
\usepackage{nicefrac}
\usepackage{latexsym}
\usepackage{mathtools}
\usepackage[noadjust]{cite}
\usepackage{paralist}
\usepackage{longtable}
\usepackage{float}
\usepackage{tikz, tikz-cd}
\usetikzlibrary{patterns}
\usepackage{hhline}
\usepackage{import}

\usepackage{enumitem}

\theoremstyle{definition}

\numberwithin{equation}{section}
\usepackage[plainpages=false,pdfpagelabels,backref=page,
citecolor=red]{hyperref}
\usepackage{xcolor}
\hypersetup{
 colorlinks,
 linkcolor={blue!80!black},
 citecolor={red!80!black},
 urlcolor={green!40!black}
 } 

\usepackage{color}

\definecolor{gr}{rgb}   {0., 0.8, 0. } 
\definecolor{bl}{rgb}   {0., 0.5, 1. } 
\definecolor{mg}{rgb}   {0.7, 0., 0.7}

\newcommand{\IC}{\mathbb{C}}

\newcommand{\R}{\mathbb{R}}

\newcommand{\cD}{\mathcal{D}} 
\newcommand{\cH}{\mathcal{H}}

\newcommand{\cL}{\mathcal{L}}

\newcommand{\cS}{\mathcal{S}} 

\renewcommand{\L}{\operatorname{L}} 
\renewcommand{\H}{\operatorname{H}} 
\newcommand{\W}{\operatorname{W}}
\DeclareRobustCommand{\Wdot}{\dot{\W}\protect{\vphantom{W}}} 
\DeclareRobustCommand{\Wdot}{\dot{\W}\protect{\vphantom{W}}} 

\newcommand{\wt}{\widetilde}

\newcommand{\dhalf}{D_t^{1/2}} 
\newcommand{\HT}{H_t} 

\newcommand{\ta}{{\scriptscriptstyle \parallel}}
\newcommand{\no}{{\scriptscriptstyle\perp}}

\newcommand{\gradx}{\nabla_x}
\renewcommand{\div}{\operatorname{div}}
\newcommand{\curl}{\operatorname{curl}}
\newcommand{\e}{\mathrm{e}} 

\renewcommand{\d}{\, \mathrm{d}} 
\renewcommand\Re{\operatorname{Re}}

\newcommand{\Le}{\mathcal{L}}

\def\Xint#1{\mathchoice
{\XXint\displaystyle\textstyle{#1}}%
{\XXint\textstyle\scriptstyle{#1}}%
{\XXint\scriptstyle\scriptscriptstyle{#1}}%
{\XXint\scriptscriptstyle%
\scriptscriptstyle{#1}}%
\!\int}
\def\XXint#1#2#3{{\setbox0=\hbox{$#1{#2#3}{%
\int}$ }
\vcenter{\hbox{$#2#3$ }}\kern-.6\wd0}}
\def\barint{\,\Xint -} 
\def\bariint{\barint_{} \kern-.4em \barint}
\def\bariiint{\bariint_{} \kern-.4em \barint}
\renewcommand{\iint}{\int_{}\kern-.34em \int} 
\renewcommand{\iiint}{\iint_{}\kern-.34em \int} 
\title[Impact of Yves Meyer's work on Kato's conjecture]{Impact of Yves Meyer's work on Kato's conjecture}
\author{Pascal Auscher}
\address{Universit\'e Paris-Saclay, CNRS, Laboratoire de Math\'{e}matiques d'Orsay, 91405 Orsay, France}
\email{pascal.auscher@universite-paris-saclay.fr}
\date{October 31, 2023}
\dedicatory{}
\allowdisplaybreaks
\begin{document}

\begin{abstract}
We discuss how the works of Yves Meyer, together with Raphy Coifman, on Calder\'on's program and singular integrals with minimal smoothness in the seventies,   paved the way not only to a solution to Kato's conjecture for square roots of elliptic operators, but also to major developments  in elliptic and  parabolic boundary value problems with rough coefficients on rough domains.  
\end{abstract}

\maketitle


\section{Introduction}

\noindent Tosio Kato's square root conjecture is one example of a question arising from one field,  formulated in a second one  and finding its solution in a third one. Namely, the question arising  in the sixties from the work of T. Kato, motivated by partial differential equations in inhomogeneous media, was set using a functional analysis framework and it was finally methods from real harmonic analysis  that put a final end to the problem as posed. 

It is likely that  R. Coifman and Y. Meyer were not aware of Kato's conjecture before the visit of A. McIntosh to Paris in 1980. McIntosh had  remarked the analogies between the multilinear series that he and they were using. McIntosh was indeed trying to attack Kato's conjecture by this tool.  
 R. Coifman and Y. Meyer were trying to prove the boundedness of  the Cauchy singular integral operator on Lipschitz curves, as part of  the program of A. Calder\'on on functional analysis and operator  algebras. This encounter allowed them to find the way of summing the series  for both the Cauchy integral and the square root in one dimension. This appears in the famous article \cite{CMcM} published in 1982. It opened the door to new developments in the above mentioned three fields and other ones, where improvements came from round trips from Cauchy integrals and Kato's conjecture, with R. Coifman, A. McIntosh and Y. Meyer as the first architects. 
  
From the point of view of Yves, it all comes from Calder\'on's insight.
Let us quote his article on the work of Calder\'on in complex analysis and operator theory \cite{Me}: 

``\textit{He [Calder\'on] pointed up a concise question which happened to be the magic key opening all doors.}''

Motivated by the example of the Cauchy integral on a curve, this concise question was to prove the $\L^2$ boundedness of what we now call anti-symmetric singular integrals. It turns out that  it is not always true and a necessary and sufficient criterion for a positive answer is the celebrated $T(1)$ theorem of G. David and J.L. Journ\'e \cite{DJ}, followed by the $T(b)$ theorem of the same authors with S. Semmes \cite{DJS}. These two theorems owe a great deal to the pioneer work of R. Coifman, Y. Meyer and A. McIntosh. Kato's conjecture one dimensional solution was also very much instrumental in the first version of the $T(b)$ theorem \cite{McM}.

Boundedness of singular integrals could be addressed in several dimension but Kato's question in higher dimension was more difficult since the direct connection with the Cauchy integral and singular integrals is lost.  Although various different proofs for the Cauchy integral's boundedness came soon after \cite{CMcM}, it took another twenty years to fully understand Kato's question in all dimension.   A new approach was needed.  We shall explain how the progress in the understanding of $T(b)$ theorems helped.

Among the applications, the solution to the square root problem could address the original questions of Kato on propagation of waves and stability with respect to the conductivity matrix of the medium. More applications connecting Kato's question and boundary value problems were also brought to light by C. Kenig \cite{Ke}. 

But the story was not finished. Almost another twenty years after the solution of  Kato's conjecture for elliptic operators, the estimates that come with and the ideas surrounding the argument found unforeseen applications in boundary value problems for elliptic and parabolic operators and Hardy space theory. It is still an active topic.

The objective of this article is to give a non-exhaustive feeling for the advances permitted by the breakthroughs made by Yves Meyer and co-authors towards partial differential equations in the last forty years with Kato's square root problem as a central theme. The kind and friendly shadow of Yves spreads all over these developments.

\section{The square root problem for elliptic operators}

\noindent  Let $A=(a_{jk})$ be a matrix of bounded measurable complex functions on $\R^n$ satisfying the accretivity inequality,
\begin{align}
\label{eq: Garding d}
\Re \sum_{j,k=1}^n a_{jk}(x) z_{k}\overline{z_{j}}  \geq \lambda |z|^2 \quad (z \in \IC^n, x\in \R^n),
\end{align}
for some $\lambda>0$.
One can then define a maximal accretive operator in $\L^2(\R^n)$ by 
\begin{align}
\label{eq:L}
 L(u)\coloneqq - \div_x (A \nabla_xu)\end{align} 
 with domain $\cD(L)=\{u\in \W^{1,2}(\R^n); \div_x (A \nabla_x u) \in \L^2(\R^n)\}$, where $\W^{1,2}(\R^n)$ denotes the usual $\L^2$-Sobolev space of order $1$. Functional analytic methods furnish fractional powers $L^\alpha$, $\alpha\in \R$, to such operators and the case of interest is $0<\alpha <1$. For $\alpha<1/2$, Kato was able to identify the domain of $L^\alpha$ as the Sobolev space $\W^{2\alpha,2}(\R^n)$ by interpolation methods. The situation when   $\alpha \ge 1/2$ remained unsolved \cite{K1}. For $\alpha =1/2$,\footnote{When $\alpha>1/2$, one does not expect a possible identification with of the domain of $L^\alpha$ with a concrete function space when $A$ is only measurable.}    Kato's conjecture (he claimed to have never explicitly  formulated  this conjecture but that it was extrapolated from his work by A. McIntosh, see \cite{Mc4}) is that the domain of $L^{1/2}$ identifies to  $\W^{1,2}(\R^n)$ with the equivalence
 \begin{align}
\label{eq:Kato}
\|L^{1/2} u\|_2 \simeq \|\nabla_x u\|_2. 
\end{align}
An interpretation in physics terms is that the kinetic energy of a wave propagating across a medium with conductivity coefficients given by $A$ is related to its spectral frequency, where spectrum is understood in the sense of the spectral theory for the operator $L$. When it is related to physics, the matrix $A$ is symmetric with real-valued entries and \eqref{eq:Kato} is a consequence of spectral theory as $L$ is self-adjoint. But Kato wanted to understand the stability of the upper inequality  in \eqref{eq:Kato} with respect to  $\L^\infty$ perturbations of $A$: that is, small perturbation in the medium should imply small perturbation in the spectrum in the sense that
\begin{align}
\label{eq:Katopert}
\|L_{0}^{1/2}u -L_{1}^{1/2} u\|_2 \lesssim \|A_{0}-A_{1}\|_{\infty} \|\nabla_x u\|_2 
\end{align} 
if $\|A_{0}-A_{1}\|_{\infty}$ is small. 
This is related to the program described in his monumental book \cite{K}. This is why the question is formulated for complex matrices:  proving  \eqref{eq:Kato} for complex entries using only $\L^\infty$ information on $A$ would automatically give analytic regularity, hence \eqref{eq:Katopert}.

The construction of $L$ follows from the representation 
\begin{align}
\label{eq:representation}
\int_{\R^n}  Lu(x) \overline{v(x)}\, \d x= \beta(u,v) \quad (u\in D(L) , v\in \W^{1,2}(\R^n))
  \end{align}
with the sesquilinear form  
  \begin{align}
\label{eq: sesquilinear form}
\beta(u,v)\coloneqq \int_{\R^n}  \sum_{j,k=1}^n a_{jk}(x) \partial_{k}u(x) \overline{\partial_{j}v(x)}\, \d x \quad (u, v\in \W^{1,2}(\R^n)).  \end{align}
This construction of the operator falls within the theory of maximally accretive sesquilinear forms. If one forgets about elliptic operators and deals with an abstract  form $B(u,v)=\langle A S u, S v\rangle$ on a Hilbert space,  where $A$ is a bounded and strictly accretive operator and $S$ a self-adjoint operator, Kato's question reformulates as follows:

{\it Is the  domain of the square root of the maximal accretive operator 
$S^*AS$ obtained from a maximally accretive sesquilinear form $B$  equal to the domain of the form, that is the domain of $S$?}

  A. McIntosh found a counter-example \cite{Mc2}, indicating that functional analytic methods would not suffice to prove \eqref{eq:Kato} in full generality. But the question  for  the differential operators in \eqref{eq:representation} coming from the forms in \eqref{eq: sesquilinear form} remained.

\section{The  approach of Coifman, McIntosh, Meyer} 

  We feel interesting to modestly plagiate the presentation by Yves in  \cite{Me}  as we want to give a stress on the Kato problem while he focused more on singular integrals. He writes

``\textit{Morever, the same key [question posed by A. Calder\'on]  led to striking discoveries in real and complex analysis.}'' 

As mentioned, this central question included the  $\L^2$ boundedness of the Cauchy integral operator. Some  other striking discoveries came later in geometric measure theory: solution to the Vitushkin conjecture by G. David \cite{Dav} and to the Painlev\'e conjecture by X. Tolsa \cite{T}. {See the article of J. Verdera in the same volume.} 

The Cauchy operator $C_{\Gamma}$ can be defined as follows: $\Gamma$ is here a Jordan Lipschitz curve in the plane,    $s\mapsto z(s)$  a parametrization and one sets, using the notation  $\mathrm{p.v.}$ for the  principal value, 
\begin{equation}\label{eq:Cauchyop}
C_{\Gamma}f(z(s))= \mathrm{p.v.}\frac{1}{\pi i}\!\int_{\R} \frac {f(z(t))z'(t)}{z(s)-z(t)} \, \d t.
\end{equation}
The approach was to expand the kernel in a multilinear series, to prove $\L^2$ boundedness for the general terms and to control the radius of convergence. Multilinear series was a topic that Coifman and Meyer were exploring at this  time using their so-called $P_{t},Q_{t}$ spectral approach and which was explained in their Ast\'erisque opus \cite{CMast}, see also \cite{CM}. When $\Gamma$ is the graph of a  Lipschitz function $A$, choosing  $z(x)=x+iA(x)$ and forgetting the harmless factor $z'(x)$, for $f$ now defined on $\R$, one studies
$$
C_{A}f(x)=  \mathrm{p.v.}\frac{1}{\pi i}\!\int_{\R} \frac {f(y)}{x-y}\left(1+i \frac{A(x)-A(y)}{x-y}\right)^{-1} \, \d y
$$
which can   be expanded as $\sum_{k=0}^\infty (-i)^k T_{k,A}f$ where  $T_{k,A}$ is the $k$th iterated Calder\'on commutator
$$
T_{k,A}f(x)=  \mathrm{p.v.}\frac{1}{\pi i}\!\int_{\R} \frac {f(y)}{x-y}\left( \frac{A(x)-A(y)}{x-y}\right)^{k} \, \d y.
$$
The term $k=0$ is the well-known Hilbert transform. A. Calder\'on was able to prove $\L^2(\R)$ operator bounds of the order $C^k\|A'\|_{\infty}^k$ for  the $k$th term with some unspecified, possibly large $C$, thus imposing $\|A'\|_{\infty}$ be small for the convergence of the series.

At the same time, A. McIntosh was exploring the possibility of proving Kato's question for differential operators.  One of the many formulas to compute the square root is
\begin{eqnarray}
\label{eqn5} L^{1/2}u = \frac{2}{\pi}
\int_0^{\infty} (1+t^2L)^{-1}Lu \, \d t, \quad u \in \mathcal{D}(L). 
\end{eqnarray} 
In one  dimension the operator takes the form $L= DaD$ where $D=-i\frac{d}{dx}$ is self-adjoint (self-adjointess turned out to be important feature of the one dimensional case with no equivalent in higher dimension) and by abuse of notation we identify the function $a$ and the multiplication by $a$, which is here an accretive function, that is a bounded complex-valued function with  real part bounded below by $\lambda>0$.
Expanding in power series with $a=1-m$, $\|m\|_{\infty}<1$, yields
$$
L^{1/2}u=  \sum_{k=0}^\infty \frac{2}{\pi} \int_0^{\infty} Q_{t}(m(1-P_{t}))^k (aDu)\,  \frac{ \d t}{t}\coloneqq \sum_{k=0}^\infty S_{k}(a	Du)
$$
with $P_{t}= (1+t^2D^2)^{-1}$ and $Q_{t}=tDP_{t}$. The strategy is to prove the boundedness of $S_{k}$ by iteration and to control the growth in $k$. Note that the first term $S_{1}$ contains $Q_{t}(m(1-P_{t}))= Q_{t}m - Q_{t} m P_{t}$. Handling the term with $Q_{t} m P_{t}$ is  already difficult and is what Coifman and Meyer could do.  Actually, A. McIntosh showed that the $k$th iterated commutator of Calder\'on can be represented as
$$T_{k,A} = \mathrm{p.v.}\frac{1}{\pi i}\!\int_{\R} R_{t}(A'R_{t})^k\,  \frac{ \d t}{t}
$$
where $R_{t}=P_{t}-iQ_{t}= (1+itD)^{-1}$.  Using this and taking into account evenness (in $t$) of $P_{t}$ and oddness of $Q_{t}$, this is in fact very close to the formula for $S_{k}$ so that the same technology can be used. For example, 
$$T_{1,A} = -\mathrm{p.v.}\frac{1}{\pi }\!\int_{\R} (Q_{t}A'P_{t}+ P_{t}A'Q_{t})\,  \frac{ \d t}{t}.
$$
This is a term that the Coifman-Meyer technique could handle. 

A.~McIntosh told me that the first proof of the $\L^2$ boundedeness of the Cauchy integral was a reduction from the solution of Kato's question for the square root of $L=DaD$, that is,  
\begin{equation}
\label{eq:rootDaD}
\|(DaD)^{1/2}u\|_{2} \simeq \|Du\|_{2}.
\end{equation}
This is not the way it is presented in \cite{CMcM}. Actually, it has been shown later that the two results are  equivalent, and we shall come back to this.  R. Coifman, A. McIntosh and Y. Meyer eventually  obtained in  the bound $C_{0}(1+k)^4\|A'\|_{\infty}^k$ for the operator norm of $T_{k,A}$ using this representation, and a renormalisation showed that the assumption $\|A'\|_{\infty}<1$ could be removed. 

This is no accident that Cauchy integrals and Kato's square root belong to the same family of operators, as evidenced by the following identity. If one considers  the operator $-i \frac{d}{dz}$ acting on smooth functions defined  on the Lipschitz curve $\Gamma$, and use parametrization $z=z(x)$, $x\in\R$, then it can be shown it is similar to $aD$ with $a(x):=1/z'(x)$,  $D=-i\frac{d}{dx}$, and a calculation shows that  
\begin{equation}\label{eq:aDaD}
 (aDaD)^{1/2} (u)(x)= \mathrm{p.v.}\frac{1}{\pi i}\!\int_{\R} \frac {Du(y)}{z(y)-z(x)} \, \d y,
\end{equation}
where the definition of the square root given by \eqref{eqn5} relies on the theory of sectorial operators that includes maximal accretive operators as a particular case.  In different notation, this integral  is the same as the one in \eqref{eq:Cauchyop}.

Pushing the techniques of \cite{CMcM} further, C. Kenig and Y. Meyer  wrote in the title of \cite{KM} that Cauchy integrals and Kato's square roots are the same (in one dimension). They proved 
\begin{equation}
\label{eq:rootbDaD}
\|(bDaD)^{1/2}u\|_{2} \simeq \|Du\|_{2}.
\end{equation} 
 by showing that the multilinear expansions  obtained for $(bDaD)^{1/2}$, with $b$ having the same properties as $a$, can be controlled with polynomial growth.   
 
  The Cauchy integral is an example of a Calder\'on-Zygmund operator. One can therefore push the ressemblance further and ask about the relation between $(bDaD)^{1/2}$ and Calder\'on-Zygmund operators. Later,  P. Tchamitchian and I showed in \cite{AT} that even when $b\ne a$ so in particular when $b=1$,  $(bDaD)^{1/2}=RD$ where $R$ is a Calder\'on-Zygmund operator. The  $\L^2$ boundedness  did not use multilinear expansions;  it is a (simple) consequence of  the $T(b)$ theorem  of G. David, J.L. Journ\'e and S. Semmes. Note that this  reproved Kato's conjecture as a particular case, and reversed the chronology of ideas. Indeed,  recall that   the first version of a $T(b)$ theorem is due to A. McIntosh and Y. Meyer \cite{McM} using ideas from the proof of the Kato's conjecture. 
   
   The ideas surrounding the $T(b)$ theorem became the core of all further developments as we shall see.

\section{The solution of the square root problem for elliptic operators}

The solution in higher dimension came in several steps that took twenty years. In his 1990
survey article \cite{Mc4}, A. McIntosh wrote 
``\textit{It remains a challenge, however, to prove the multilinear estimates needed to solve the square root problem of Kato for elliptic operators, or to find an alternative approach.}''

First, one had to give up the idea that the multilinear series that worked in one dimension would do the same in higher dimension. Such techniques were worked out by R. Coifman, D. Deng, Y. Meyer \cite{CDM} and independently by E. Fabes, D. Jerison and C. Kenig \cite{FJK1}. They were optimized by J. L. Journ\'e \cite{J}. The best bound he obtained is that, writing $A^{-1}=I-M$ with $\| |M| \|_{\infty}<1$ (the $\L^\infty$ norm of the matrix $M(x)$ as an operator on $\IC^n$ with hermitian structure),  the $k$th term in the multilinear expansion  is controlled by  $C^k\| |M| \|_{\infty}^k$ with $C$ on the order of $n^{1/2}$. This is not enough.

At the same time, while the Cauchy integral $\L^2$ boundedness was given many new proofs essentially because of the tight relation with holomorphic functions, there were no new ones for the one-dimensional Kato problem until the nineties. One was using an adapted wavelet basis \cite{AT1} but the same idea in higher dimension remains mysterious.   The proof in \cite{AT} around the mid nineties showed that there was hope for a proof without doing such expansions by incorporating the $T(b)$ technology.  But this proof (in one dimension) used the full force of Calder\'on-Zygmund theory.  An example of C. Kenig, see \cite{AT2}, showed  that  square roots  in higher dimension cannot be related to  Calder\'on-Zygmund operators as some  $\L^p$ estimates usually obtained as a consequence of the Calder\'on-Zygmund extrapolation method fail. 

Nevertheless, in any dimension, S. Semmes \cite{S} proposed  a variant of the $T(b)$ theorem requiring less information than the one needed for Calder\'on-Zygmund theory.   Let us explain this and, at the same time,   describe the role of $b$. This variant was a criterion to prove square function estimates of the form
\begin{align}
\label{eq:thetat}
\int_{0}^\infty \|\theta_{t}f\|_{2}^2\, \frac{\d t}{t} \lesssim \|f\|_{2}^2
\end{align} when $(\theta_{t})_{t>0}$ is a family of operators acting on complex functions, with kernels $\theta_{t}(x,y)$ having good pointwise   bounds and some  regularity  in the $y$ variable.
In that case, it can be shown that \eqref{eq:thetat} holds if and only if $ |\theta_{t}1(x)|^2 \frac{\d x \d t}{t}$ is a Carleson measure, which means that there is a constant $C<\infty$ such that for all cubes $Q\subset \R^n$ with sides parallel to the axes (Euclidean balls work as well), $\ell(Q)$ being its sidelength and $|Q|$ its volume,
\begin{equation}\label{eq:Carl}
\int_Q \!
\int_0^{\ell(Q)}| \theta_t1(x)|^2 \,
\frac{\d x\d t}{t}
\le C |Q|.
\end{equation}
Here, $\theta_{t}1$ is the operator $\theta_{t}$ applied to the constant function 1.  But the calculation of $\theta_{t}1$ might just be impossible so this could be useless. However, this equivalence contains the fact that the inequality
$$
\int_0^\infty\!\int_{\R^n} |\theta_t f(x)-\theta_t1(x)\cdot
S_tf(x)|^2
\frac{\d x\d t}{t}
\le C \|f\|_{2}^2
$$
always holds,
where $S_{t}$ is a dyadic martingale ($S_{t}f(x)$ is the average of the function $f$ on the dyadic cube containing $x$ with  $\ell(Q)/2\le t <\ell(Q)$). In passing, this reduction to a principal part $\theta_t1(x)\cdot
S_tf(x)$ is done using the $P_{t},Q_{t}$ techniques of Coifman and Meyer. So if one knows functions $b$ for which $|\theta_tb(x)|^2 \,
\frac{\d x\d t}{t}$ is a Carleson measure  and the real part of $S_{t}b(x)$ is bounded below, at least locally on each Carleson window $Q\times (0,\ell(Q))$,  one can conclude that \eqref{eq:Carl} holds.  In the one dimensional application towards Kato's problem, where we recall that  $a$ is  a bounded complex-valued function with  real part bounded below by a positive number, one has  $\theta_{t}f= (1+t^2L)^{-1} (taf)'$ so that  $ \theta_{t}a^{-1}=0$  and  the above algorithm works out perfectly. 

The approach of Semmes was inspiring for square roots in higher dimension because some easy reductions from  functional analysis tell us that the upper estimate in
\eqref{eq:Kato} suffices and that this one is equivalent to \eqref{eq:thetat} with $\theta_{t}F= (1+t^2L)^{-1} t\div_{x}(AF)$, $F$ being here $\IC^n$-valued,  namely
\begin{align} \label{eq:SFEKato}
 \int_0^\infty\!\int_{\R^n} |(1+t^2L)^{-1} t\div_{x}(AF)(x)|^2
\frac{\d x\d t}{t}
\le C \|F\|_{2}^2.
\end{align}
Despite the fact that formally $\theta_{t}A^{-1}(x)=0$,  one was missing pointwise bounds, $y$-regularity in any sense, and   the facts that $\theta_{t}I$ is $\IC^n$-valued ($I$ being the identity matrix) and that $A$ is  a matrix complicated the algebra. Many reasons to stop there! 

Nevertheless,  a possible strategy was developed  in the Ast\'erisque book with P. Tchamitchian \cite{AT2}. Basically we proved in this rough context   (under some further technical hypotheses) that the same reduction to \eqref{eq:Carl} as in Semmes' work holds together with the valid estimate
$$
\int_0^\infty\!\int_{\R^n} |\theta_t \nabla_{x} f(x)- \theta_tI(x)\cdot
S_t\nabla_{x} f(x)|^2
\frac{\d x\d t}{t}
\le C \|\nabla_{x} f\|_{2}^2.
$$
Here, writing $\theta_tI(x)= (\gamma_{t,1}(x), \ldots, \gamma_{t,n}(x))$,  $$\theta_tI(x)\cdot
S_t\nabla_{x} f(x)=\gamma_{t,1}(x)
S_t(\partial_{x_{1}} f)(x)+\cdots+ \gamma_{t,n}(x)
S_t(\partial_{x_{n}} f)(x).$$
This led us to formulate a $T(b)$ criterion for square roots with the new idea that the test function to construct could be not just one function but a family of functions $(F_{Q})$ indexed by cubes: one would need   \begin{itemize}
\item an $\L^2$ control on their gradients on an enlargement of the cube $Q$,
\item that $F_{Q}$ is adapted to $L$, in the sense that $|\theta_{t}\nabla_{x}F_{Q}(x)|^2\frac{\d x\d t}{t}$ is controlled on an enlargement of  the Carleson window $Q\times (0,\ell(Q))$,
\item and that   $\theta_tI(x)\cdot
S_t\nabla_{x} F_{Q}(x)$ controls $|\theta_tI(x)|$ on  the Carleson window.
\end{itemize}
 The first two requirements are not so difficult to achieve but the last one is complicated because $ \theta_tI(x)$ is $\IC^n$-valued, and also 
because the control is not on  the Carleson window in general but on a  subregion. Fortunately, the size of missing region where there is no control can be estimated using  stopping-time arguments (this is where  dyadic martingales are very useful).  

This $T(b)$ criterion for square roots turned out to be the successful approach. S. Hofmann and A. McIntosh \cite{HMc} were able to apply this criterion to prove Kato's conjecture in dimension $n=2$. The  idea  brought by M. Lacey for arbitrary dimension  in \cite{HLMc} was to use a sectorial decomposition of $\IC^n$ to force  $\theta_tI(x)$ to take values in small cones about a finite number of unit vectors $w$ by cut-off and, for each fixed $w$, to choose the $F_{Q}$ so that $\nabla F_{Q}$ is close to   the  direction  $\overline w$ in average on the cube $Q$. Removal of the assumption of pointwise bounds on the kernels of the heat semigroup $\e^{-t^2L}$ in \cite{HLMc} was done in the final assault \cite{AHLMcT}. 
 
Let me mention that M. Christ  had  already the idea to use varying test functions with $\L^\infty$ control to prove a ``local'' $T(b)$ theorem for singular integrals in his article about analytic capacity \cite{Ch}. His proof does not allow for a kind of $\L^p$ control with $p$ finite though. The possibility of having $\L^p$ control with finite $p$ generated some search to refine the exponents and testing conditions in local $T(b)$ theory for singular integrals, motivated by this use in partial differential equations. The best (and nearly optimal) result is due to T. Hyt\"onen and F. Nazarov \cite{HN}.

The strategy of proof of  Kato's square root problem is very flexible and allows generalizations where the proof adapts:  elliptic systems of arbitrary high order  possibly having lower order terms with bounded measurable coefficients (or even unbounded in some natural spaces in order to preserve ellipticity);  Ellipticity can be relaxed to a G\aa rding inequality (which is weaker than the pointwise ellipticity \eqref{eq: Garding d}). Lower order terms do  bring  technical difficulties that are minor compared to the treatment of the highest order terms. Ellipticity can also degenerate in weighted sense, see \cite{CUR}.

\section{The first order approach}

Even before  Kato's question was solved in higher dimension, a first order approach was explored jointly with A. McIntosh and A. Nahmod in the case of the one dimensional setting \cite{AMcN1}. The idea is to see the second order operator 
$L= -b(x)\frac{d}{dx} (a(x) \frac{d}{dx})$ (the same as $bDaD$ above) as (part of) the square $T^2$  of a first order differential system $T$ whose domain is $\W^{1,2}(\R; \IC^2)$ and whose spectrum belongs to a bisector $|\arg(\pm z)|\le \omega_{T} <\pi/2$ (such operators belong to the class of bisectorial operators). 
 Then identification of  the domain of $L^{1/2}$ follows from knowing the one of  $(T^2)^{1/2}$ and one hopes for
  \begin{align}
 \label{eq:boundedness}
 \|(T^2)^{1/2} u \|_{2}\sim \|Tu\|_{2}
 \end{align}
for all $u$ in the domain of $T$.  It amounts to showing that $(T^2)^{1/2}T^{-1}$, restricted to the range of $T$,  is bounded on $\L^2$.
   The advantage of this point of view is that all can be seen within the holomorphic functional calculus of $T$ and for this A. McIntosh had established a criterion using quadratic estimates (see \cite{Mc}).  More precisely, he established that bounded and holomorphic functions of bisectorial operators  (this was done for sectorial operators but the proof is the same) $T$ on larger bisectors than the one given by $\omega_{T}$, generate bounded operators if and only if there is equivalence between a certain square function norm and the ambiant norm:
 \begin{align}
 \label{eq:SFE}
  \int_0^\infty\| {\lambda{}}T(1+{\lambda{}^2}T^2)^{-1} h \|_{2}^2 \, \frac{\mathrm{d}\lambda{}}\lambda{} \simeq \|h\|_{2}^2
\end{align}
for all $h$ in the closure of the range of $T$.  The function  $(z^2)^{1/2}z^{-1}$ being bounded and holomorphic, this shows that  \eqref{eq:boundedness}  follows from \eqref{eq:SFE}. 
Here, $T={\bf BD}$ with 
  \begin{align*} {\bf B}= \begin{bmatrix} b & 0 \\ 0 & a \end{bmatrix}
, \quad 
{\bf D}= \begin{bmatrix} 0&\frac{d}{dx} \\ -\frac{d}{dx} & 0 \end{bmatrix}
\end{align*} so that 
$$T^2=  \begin{bmatrix}  -b(x)\frac{d}{dx} (a(x) \frac{d}{dx}) & 0 \\ 0 &  -a(x)\frac{d}{dx} (b(x) \frac{d}{dx}) \end{bmatrix}=\begin{bmatrix}  bD aD & 0 \\ 0 &  aDbD \end{bmatrix} \ . $$

In one dimension, one can obtain kernel bounds for the operators ${\lambda{}}T(1+{\lambda{}^2}T^2)^{-1}$ so that the square function estimate \eqref{eq:SFE} can be shown following  the Semmes version of  $T(b)$ theorem for square functions. 
 
We pause to come back to the affirmation that Kato's square root and Cauchy integrals are the same. Using this setup one can indeed show  without computing integral formulas as before, but only using  functional analytic tools such as interpolation,   that \eqref{eq:rootDaD} for all $a$,    \eqref{eq:rootbDaD} for all $a,b$, and  \eqref{eq:rootbDaD} for all $a=b$ are equivalent,
 see \cite{AMcN1}. Thanks to \eqref{eq:aDaD},  this  is equivalent to the $\L^2$  boundedness of the Cauchy integral operator.

 As we saw, the solution of  Kato's square root problem in all dimension follows from a $T(b)$ strategy even though one does not have pointwise kernel bounds.  A. Axelsson, S. Keith and A. McIntosh found a way to extend the first order approach presented above in one dimension to  higher dimension using analogy with metric perturbations of Dirac operators and were able to formalise the $T(b)$ argument within this framework \cite{AKMc1}. In differential geometry, the Dirac operator is $\Pi=d+d^*$ where $d$ is the exterior derivative on forms and $d^*$ its adjoint with respect to the $\L^2$ inner product. What really matters here is that $d$ is a nilpotent first order constant coefficient differential operator  and that $\Pi$ is coercive on its range: $\|\nabla u\|_{2} \lesssim \|\Pi u\|_{2}$ for $u$ in the $\L^2$ range of $\Pi$. 
Metric perturbations take the form $\Pi_{B}=d+B_{1}d^*B_{2}$ where $B_{1}, B_{2}$ are bounded operators and $B_{1}$ is the inverse of $B_{2}$ (being usually a multiplication operator). In general, it suffices to assume that $B_{1},B_{2}$ enjoy some partial coercivity (a lower bound for the real parts of $\langle B_{i}h,h\rangle$ for $h$ in appropriate subspaces)   and that $d^*B_{2}B_{1}d^*=0$ and $d B_{1}B_{2}d=0$. For applications to  Kato's problem, it suffices to set this framework on 0 and 1 forms, but further applications to solving partial differential equations for forms may require to use the full exterior algebra. 
They establish a (sufficient) criterion to prove \eqref{eq:SFE} for $T=\Pi_{B}$ when it is additionally required that $B_{1},B_{2}$ are multiplications by bounded matrix-valued functions. According to A. McIntosh's result, \eqref{eq:SFE} for $\Pi_{B}$ yields the $\L^2$ boundedness of $(\Pi_{B}^2)^{1/2}\Pi_{B}^{-1}$ defined on the range of $\Pi_{B}$.  This is a stronger statement as it implies as corollaries 
Kato's square root  conjecture for equations and systems in all dimension, the $\L^2$ boundedness of the Cauchy integral operator and its Kenig-Meyer extension \eqref{eq:rootbDaD} for all $a,b$ and more. In this framework, $d$ can be replaced by a nilpotent constant coefficient first order differential operator provided that the triple $(\Gamma, B_{1},B_{2})$ satisfies the above requirements.

\section{The Lions problem}

Jacques-Louis Lions asked whether  Kato's square root problem can be addressed on a domain for operators with mixed boundary conditions (possibly allowing lower order terms with bounded coefficients) \cite{Li}. That is, on a part of the boundary, the operator comes with Dirichlet condition and on the ``complement'' the operator comes with Neumann condition. The construction of the operator uses the theory of quadratic forms. The first positive results for Lipschitz domains and pure Dirichlet or Neumann conditions are in a work with P. Tchamitchian, roughly using localization and extension techniques to reduce to the case of $\R^n$ \cite{ATLip}.  For mixed  conditions, such techniques do not apply and it was solved by A. Axelsson, S. Keith and A. McIntosh on Lipschitz domains by implementing their first order approach \cite{AKMc2}. 

The state of the art is a result by S. Betchel, M. Egert and R. Haller-Dintelmann  in geometries beyond the Lipschitz domains:  Lions question has a positive answer on an open and possibly unbounded set in $\R^n$ under two simple geometric conditions: The Dirichlet boundary part is Ahlfors--David regular and a quantitative connectivity property in the spirit of locally uniform domains holds near the Neumann boundary part.
 For example, the Neumann boundary condition can be treated on  $(\varepsilon,\delta)$-domains, in particular on the interior of the von Koch snowflake. We refer to \cite{BEH}.

\section{Elliptic boundary value problems}

As mentioned earlier, the estimate \eqref{eq:Kato} for $L= - \div_x (A(x) \nabla_x)$ has consequences for the boundary value problems of the elliptic equation 
\begin{align}
\label{eq:block}
\partial_{t}^2u(t,x)+  \div_x (A(x) \nabla_xu(t,x))=0
\end{align} 
where $x\in \R^n$ and $t>0$ and the boundary is identified with $\R^n$. Indeed, if $f$ is a regular Dirichlet data with $\|\nabla_x f\|_{2}<\infty$ (this is why this is tagged as the regularity problem), then the formal  solution $u(t,x)=\e^{-tL^{1/2}}f(x)$, satisfies $\sup_{t>0}\|\nabla_{t,x}u\|_{2}\sim \|\nabla_x f\|_{2}$; if $g$ is a Neumann data with $\|g\|_{2}<\infty$ then the formal solution $u(t,x)=-\e^{-tL^{1/2}}(L^{-1/2}g)(x)$, has conormal derivative $g$ at $t=0$, and satisfies $\sup_{t>0}\|\nabla_{t,x}u\|_{2}\sim \|g\|_{2}$. Actually, much more can be said in terms of  non-tangential maximal functions estimates for such solutions and also that they are unique with such requirements. The state of the art, together with extension to other boundary spaces which require the theory of Hardy spaces adapted to $L$, appears in a recent work with M. Egert \cite{AE}.

Equation \eqref{eq:block} is a particular occurence, called block form, of the more general elliptic equations
\begin{align}
\label{eq:non-block}
 - \div_{t,x} (A(x) \nabla_{t,x}u(t,x))=0,
\end{align} 
where $A$ is a full $(n+1)\times (n+1)$ matrix of bounded functions independent of the transverse direction $t$ to the boundary. The ellipticity is given by \eqref{eq: Garding d} in $n+1$ dimension.  

The motivating example is a simple pullback of the Laplace equation $\Delta u=0$ above a Lipschitz graph to flatten the boundary: the solvability of the Dirichlet problem for data in $\L^2$ in such a situation was proved by B. Dahlberg \cite{Da} and the regularity and Neumann problems in $\L^2$ by D. Jerison and C. Kenig \cite{JK2} based on a Rellich identity that became central in the topic. The $\L^2$ boundedness   of the double layer potential proved by R. Coifman, A. McIntosh and Y. Meyer as a consequence of the method of rotations from the boundedness of the Cauchy integral operator could then be used, as G. Verchota showed in 1984 \cite{V}, to reprove such results by the method of layer potentials. For  real and symmetric operators of the form \eqref{eq:non-block}, solvability of the Dirichlet  problem in $\L^2$  was obtained by D. Jerison and C. Kenig \cite{JK1}, and solvability of the regularity and Neumann problems in $\L^2$ by C. Kenig and J. Pipher in the mid nineties \cite{KP}, while the method of layer potentials was developed a decade ago \cite{AAAHK} thanks to the construction of fundamental solutions by S. Hofmann and S. Kim \cite{HK}. This allowed additional perturbation results, while the Rellich identities are very rigid and do not allow for perturbations.  

In the generic situation of  \eqref{eq:non-block}, one could ask what happens  when the coefficients also depend on   $t$ but counter-examples by A. Caffarelli, E. Fabes and C. Kenig \cite{CFK} show that there is no hope for positive results without perturbative assumptions with respect to  $t$. We refer to the book by C. Kenig \cite{Ke} for the state of the art until the mid nineties but the topic has grown since.   

It was not clear at all whether the connexion between square roots and boundary value problems extends to non-block form equations. In fact, the first order approach appears to be more adapted, in the spirit of  the Cauchy-Riemann system  in two dimension and the $\div-\curl$ systems of E. Stein and G. Weiss in higher dimension for the Laplace equation \cite{steinweiss}.  

The work with A. Axelsson and S. Hofmann brought the first stone \cite{AAH}: 
the strategy in that paper begins by writing \eqref{eq:non-block} in the equivalent form
\begin{align} \label{eq:Laplacein1order}
\begin{cases}   
  \div_{t,x} A(x) F(t,x)  =0, \\
  \curl_{t,x} F(t,x) =0,
\end{cases}
\end{align}
with $F(t,x)=\nabla_{t,x}u(t,x)$ as unknown. Note that the second equation is a dummy one if one knows $u$; however, this allows one to recover a gradient among all $F$ solving the first equation. Working out the system reveals the evolution equation $\partial_{t}F+T_{A}F=0$, where $T_{A}$ is a first order differential operator related to the $\Pi_{B}$ operator as in \cite{AKMc1}. Solving the boundary value problems can be done in two steps. First solve the initial value problem for the evolution equation when $t>0$. This can be done only for initial data $F_{0}$ in a spectral subspace $\cH^+_{T_{A}}$, called Hardy space adapted to $T_{A}$: $F$ is then obtained by applying an abstract Cauchy operator for $T_{A}$ given by the semigroup evolution for $-[T_{A}]=-(T_{A}^2)^{1/2}$ and the projection $\chi_{+}(T_{A})$ onto that subspace:
$$
F=\e^{-t[T_{A}]}\chi_{+}(T_{A})F_{0}.
$$
The second step is to construct a Hardy data $F_{0}$ from the knowledge of the (Dirichlet or Neumann) boundary data: this amounts to inverting  operators from the Hardy space into the boundary spaces. Roughly, this is an abstract way of thinking the Cauchy extension from a Lipschitz graph $\Gamma$: the Cauchy integral of an $\L^2(\Gamma)$ function furnishes a holomorphic function on the complement of the graph. Its traces on the graph from above and below furnish two functions in subspaces of $\L^2$, called the (holomophic) Hardy space $\cH^\pm(\Gamma)$ by analogy with what happens in the flat case of the real line.  The $\L^2$ boundedness of the Cauchy integral operator is precisely what makes the splitting  $\L^2(\Gamma)=\cH^+(\Gamma)+ \cH^-(\Gamma)$  topological.

Each step applies  only for $A$'s being a $\L^\infty$ perturbation of  well understood cases such as constant coefficients or real symmetric (or even complex hermitian). This is because $T_{A}$ contains extra multiplication terms compared to  $\Pi_{B}$ 	and this complicates the analysis. 
 
An important simplification made with A. Axelsson and A. McIntosh around 2010 in \cite{AAMc} is that there is a better choice of  unknown function than the gradient of $u$ to do the first order approach: it is the conormal gradient $\wt F\coloneqq\nabla\!_{A}u$. The latter is a vector  function obtained by replacing in $F$ the time derivative  $\partial_{t}u$ by the conormal derivative $e_{0}\cdot A\nabla_{t,x}u$ where $e_{0}$ is the unit vector normal and upward to the boundary. The tangential part of $\wt F$ is still of gradient form (or is tangential curl free). Then \eqref{eq:non-block} is equivalent to an evolution equation
$\partial_{t}\wt F+DB\wt F=0$ where 
\begin{align}
\label{eq:DB}
D= \begin{bmatrix} 0& \div_x \\ -\nabla_x & 0 \end{bmatrix}, \quad 
B= \begin{bmatrix} I & 0 \\ A_{\ta\no} & A_{\ta\ta} \end{bmatrix}
\begin{bmatrix} A_{\no\no} & A_{\no\ta} \\ 0 & I \end{bmatrix}^{-1}
\end{align}
when writing $A= \begin{bmatrix}  A_{\no\no} & A_{\no\ta} \\ A_{\ta\no} & A_{\ta\ta} \end{bmatrix}$. The representation as  $2\times 2$ matrices  corresponds to the splitting $\IC^{1+n}=\IC\times \IC^n$ in (scalar) normal  and (vector) tangential coordinates and where $\wt F(t,\cdot)$ is required to belong to the closure of the range of $D$ for all $t>0$. The operator $D$ is  self-adjoint on $\L^2(\R^n;\IC^{1+n})$ and  $A\mapsto B$ is a self-inverse map on bounded and accretive matrices, so that $DB$ is again defined as bisectorial operator. The square function estimate for $DB$ can be proved directly once again by a $T(b)$ argument or deduced from the $\Pi_{B}$ formalism of \cite{AKMc1} (actually, it is proved in \cite{AAMc2} by  functional analytic arguments that they are equivalent up to isomorphisms,) and the gain of this $DB$ formalism over the $T_{A}$ formalism is that the first step of the strategy described above can be shown for all complex matrices $A$ (in retrospect, it does too for $T_{A}$ but it could not be seen directly); it is only for the second step that things can be more subtle and, so far, there is no full understanding of it. 

In the last ten years, the $DB$ formalism has been explored with great success in its relation to elliptic boundary value problems for equations or systems with new results on  extrapolation to other boundary spaces \cite{AusSta}, solvability results  and uniqueness \cite{AM,AE-APDE}, solvability and uniqueness for data with fractional regularity \cite{AA},  perturbation results with $t$-dependent coefficients \cite{AA1}, etc. Compared to the  technology developed for second order equations with real coefficients,  an advantage is that it completely avoids use of maximum principle to establish the needed square functions. It also imposed to streamline the overall intrication of various estimates in the field.

\section{Link between first order approach and layer potentials}

We wish to isolate a spectacular success of the first order $DB$ approach, due to A. Ros\'en \cite{R}. 

As said earlier, there is a construction of fundamental solutions for the equations \eqref{eq:non-block}  in dimension $n+1\ge 3$ assuming that solutions of such equations enjoy local boundedness and H\"older regularity in the spirit of de Giorgi's results for real equations.
Basically, the fundamental solution should be a function $\Gamma(t,x,s,y)$ with 
\begin{align*} -\div_{t,x} (A(x) \nabla_{t,x}\Gamma(t,x,s,y))&=\delta(s,y), \\
- \div_{s,y} (A^*(y) \nabla_{s,y} {\Gamma^*(s,y,t,x)})&= \delta(t,x),
 \end{align*}
where $\Gamma^*(s,y,t,x)=\overline {\Gamma(t,x,s,y)}$ and $\delta(X)$ denotes the Dirac mass at the point $X$. It is constructed in such a way to obtain bounds at $\infty$. The
 $t$-independence of the coefficients, which we assume from now on,  further implies that $\Gamma(t,x,s,y)=\Gamma(t-s,x,0, y)$.

This can be used to write out the single and double layer potentials by formally setting 
$$
S_{t}f(x)= \int_{\R^n} \Gamma(t,x,0, y)f(y)\, \d y  \quad (t\ne 0)
$$
$$D_{t}f(x)=\int_{\R^n} f(y) \ e_{0}\cdot\overline {{A^*(y)}{\nabla_{s,y}\Gamma^*(s,y,t,x)}|_{s=0}}\, \d y \quad (t>0).
$$
These expressions naturally arise for smooth coefficients when invoking the Green's formula.  In his 2006 address to the ICM \cite{H2}, S.~Hofmann formulated the conjecture that $\nabla_{t,x}S_{t}$ and $D_{t}$ should be uniformly bounded as operators on $\L^2(\R^n)$ when $t>0$ and have strong limits at $t\to 0$ with bounds only involving dimension and the ellipticity constants. It was known that such bounds hold provided the Dirichlet problem for equations \eqref{eq:non-block} is solvable  for $\L^2$ data (thus, the bounds also depend on this information). 
The question was to eliminate this \textit{a priori} information since it is not known in full generality.

This conjecture was first proved by A.~Ros\'en around 2010. He showed that  $S_{t}$ and $ D_{t}$ can be respectively identified to operators $\cS_{t}$ and $ \cD_{t}$ arising from the holomorphic functional calculus of $DB$ and   $BD$ respectively, with $D$ and $B$ defined in \eqref{eq:DB}. Consequently, he showed that $\nabla_{t,x}\cS_{t}$ and $\cD_{t}$ are bounded for \textbf{all} complex   $A(x)$ as in \eqref{eq:non-block}.  He is neither assuming  local boundedness nor  H\"older regularity of solutions. It is only  when these assumptions are made that the single and double layer operators have an integral representation as above. 

Let us explain Ros\'en's idea. Writing $\cL u= - \div_{t,x} (A(x) \nabla_{t,x}u) $ as an operator from the homogeneous Sobolev space $\Wdot^{1,2}(\R^{n+1})$ into its dual, it is invertible and the fundamental solution can  somehow be interpreted as a distributional kernel of $\cL^{-1}$. But $\cL^{-1}$ and the single layer operator $\cS_{t}$  of Ros\'en are  related as follows. 
Assume $f\in C_{0}^\infty(\R^{n+1})$. If $n=1,2$ assume furthermore that $f=\div_{x}F$ for some $F\in C_{0}^\infty(\R^{n+1}; \IC^{n})$. Then, he proves 
\begin{align}
\label{eq:cvgamma1}
(\cL ^{-1}f)(t,x)
= \mathrm{p.v.}\!\int_{\R} \cS_{t-s}  f_s(x) \ ds, \quad \mathrm{in }\  \dot W^{1,2}(\R^{n+1}),
\end{align}
 where $f_{s}(x)=f(s,x)$. Thus, one can construct   $\cL^{-1}$ (and prove estimates for the fundamental solution)  from knowledge of $\cS_{t}$.   Conversely,  $\Le ^{-1}$ completely determines $\cS_{t}$ (recall we assume $t$-independence here) exactly as the fundamental solution determines the single layer potential by fixing $s=0$.   More precisely, let  $f\in C_{0}^\infty(\R^{n})$ with $\int_{\R^n} f(x)\, \d x=0$ if $n=1,2$ and  $\chi_{\varepsilon}(s)= \frac{1}{\varepsilon} \chi(\frac{s}{\varepsilon})$ with $\varepsilon>0$ and $\chi\in C_{0}^\infty(\R)$ satisfying $\int_{\R}\chi(s)\, ds=1$. Set $f_{\varepsilon}(s,x) := \chi_{\varepsilon}(s)f(x)$. Then  for all $t\in \R$, $(\Le ^{-1}f_{\varepsilon})_{t}$ converges in $\dot W^{1,2}(\R^n)$ to $\cS_{t}  f$ as $\varepsilon\to 0$ \cite{AE-APDE}. Hence,  knowledge of $\Le ^{-1}$ is dictated by the $DB$-functional calculus, and the approach of Ros\'en is natural once one knows that this holomorphic functional calculus is bounded! 

\section{The  Dirichlet problem for real but non-symmetric elliptic operators}

As said above,  the  $\L^2$ solvability of the Dirichlet problem   for real and symmetric elliptic operators  as in \eqref{eq:non-block}  had been known since the early eighties in \cite{JK1}. The symmetry assumption was used to prove by a clever argument and integration by parts an integral identity which implies control of  the conormal derivative by the tangential derivative in $\L^2(d\sigma)$ norm on the boundary. By known manipulations, this identity yields the desired estimates.

C. Kenig \cite{Ke} asked whether one can  drop the symmetry assumption. But then, the integral identity is no longer available. At least, can one prove that $L$-harmonic measure  is  $A_{\infty}$ with respect to Lebesgue measure   on the boundary?   This fact implies  $\L^p$ solvability of the Dirichlet problem for some large unspecified exponent $p$. The $L$-harmonic measure is not just one measure. It is in fact the family of Borel measures $\omega_{X}$, indexed by points $X$ of the upper-half space,  obtained by applying the Riesz representation theorem to  the maps sending any bounded non-negative compactly supported continuous function $f$ on the boundary to the value at $X$ of the (unique)  non-negative solution with data $f$. These measures are all mutually absolutely continuous thanks to the Harnack principle, but this is not (yet) enough to conclude for the $A^\infty$ property with respect to surface (Lebesgue) measure.   

With H. Koch, J. Pipher and T. Toro \cite{KKPT}, he was able to give  a positive answer in dimension $n+1=2$ in 2000, transposing  ideas on $\varepsilon$-approximability of harmonic functions to solutions of the equation and using transformations of the operator that can only be done in two dimension.

In 2014, C. Kenig, S. Hofmann, S. Mayboroda and J. Pipher   settled this question in all dimension \cite{HKMP}. Let us quote their article: ``The non-symmetric case, however, was not achievable via previously devised methods (\ldots), while being important for several independent reasons. First, the well-posed results for equations with real non-symmetric coefficients and associated estimates on solutions provide the first step towards understanding the operators with complex coefficients, in the non-Hermitian case. The latter is absolutely necessary to establish analyticity of the solution as a function of the coefficients even when the coefficients are real, one of Kato's long-time goals, which would then allow one to study also hyperbolic problems (see \cite{AHLT} for the ``block matrix'' case) [see also \cite{Mc1}]. From the analytic point of view, this can be viewed as a far-reaching extension of the Kato square root problem and Kato's program. Furthermore, the equations with complex coefficients offer the simplest model of elliptic systems retaining their major difficulties, with multiple entry points to the theory of elasticity and other applications.''

Here,  the main new ingredients are the square functions related to  Kato's square root problem for $-\div_{x}A_{\ta\ta}(x) \gradx$, where $A_{\ta\ta}(x)$ appears in the $2\times 2$ matrix form as in Section 8.   The heuristic and vague idea behind the use of such estimates is  roughly  to treat \eqref{eq:non-block} in  perturbation form
\begin{align*}
 - \div_{x} (&A_{\ta\ta}(x) \nabla_{x}u(t,x))= \div_{x} (A_{\ta\no}(x) \partial_{t}u(t,x)) \\
 & +\partial_{t} (A_{\no\ta}(x) \nabla_{x}u(t,x)) + \partial_{t} (A_{\no\no}(x) \partial_{t}u(t,x)).
\end{align*} 
However, there is no smallness attached to such a decomposition.
The term  with $A_{\ta\no}$ is handled thanks to the square function bounds \eqref{eq:thetat} for  $\theta_{t}=(1-t^2\div_{x} (A_{\ta\ta}(x) \nabla_{x}))^{-1}t\div_{x}$, while the term with $A_{\no\ta}$ requires these bounds  with $\theta_{t}=(1-t^2\div_{x} (A_{\ta\ta}^*(x) \nabla_{x}))^{-1}t\div_{x}$. The last term offers no special difficulty. 

\section{The square root problem for parabolic operators}

Parabolic equations and elliptic equations usually are closely related. For example, if one considers 
\begin{align}
\label{eq:H}
 H(u)\coloneqq  \partial_{t}u- \div_x (A(x,t) \nabla_xu)
 \end{align} 
 for $t\in \R$ and $x\in \R^n$, with $A$ bounded and satisfying \eqref{eq: Garding d} also uniformly with respect to $t$, one may ask whether the technology seen above for elliptic operators $L$  transfers to parabolic operators  $H$. Until recently, no results were available when the $t$-regularity of the coefficients is merely measurable and one often has to assume at least a control on  half-order derivatives with respect to $t$. In a series of papers, K. Nystr\"om and collaborators explored such situations in order to construct tools for solving boundary value problems for parabolic equations \cite{N1, N2, CNS}. It is assumed time-independence and they show that the square function estimates   extend in parabolic setting, i.e., replacing $L$ by $H$ in \eqref{eq:SFEKato} (and the $t$ variable integration should be renamed to avoid confusion with time), using again a $T(b)$ technology.    These works opened new doors.

The challenge was to understand whether these square function estimates hold without the time-independent assumption.  In fact, these estimates are akin to show a Kato's square root estimate for $H$.\footnote{E.~Ouhabaz (private communication) told us afterwards that the time-independence allows one to  use rather well-known techniques of functional analysis involving interpolation and maximal regularity to obtain \eqref{eq:parabolicKato} below. But this no longer applies in the time-independent case, and it was again fortunate that one did not know this abstract argument before.} What does that mean?

In order to set the context,  define the energy space $$\mathsf{V} := \H^{1/2}(\R; \L^2(\R^n)) \cap \L^2(\R; \W^{1,2}(\R^n)),$$ contained in $\L^2(\R^{n+1})$ with norm
$$
\|u\|_{V}= \left(\|u\|_{2}^2+\|\gradx u\|_{2}^2+ \|\dhalf u \|_{2}^2\right)^{1/2}
$$
 where $\dhalf$ is the half-order derivative in the $t$ variable defined through the Fourier transform with the multiplier $|\tau|^{1/2}$.   Then define the parabolic operator $H$ as an operator $\mathsf{V} \to \mathsf{V}^*$ via a sesquilinear form,
\begin{align*}
 \langle H u, v \rangle := \iint_{\R^{n+1}} A \nabla_x u\cdot\overline{\nabla_x v}+ \HT\dhalf u\cdot \overline{\dhalf v} \d x\d t \qquad (u,v \in \mathsf{V}),
\end{align*}
where $\HT$ is  the Hilbert transform in the $t$ variable defined in such a way that $\partial_{t}= \dhalf  \HT\dhalf$. 
Carving out the surprising analogy with elliptic operators even further, one can show that $H$ with maximal domain $ \{u \in \mathsf{V} : H u \in \L^2(\R^{n+1}) \}$ in $\L^2(\R^{n+1})$ is a maximal-accretive operator, hence it possesses a square root. With M. Egert and K. Nystr\"om \cite{AEN1}, we showed that 
 the domain of $H^{1/2}$ is that of the accretive form, that is $ \mathsf{V}$,  with the two-sided estimate
\begin{align}\label{eq:parabolicKato}
\|H^{1/2}\, u\|_{2} \sim \|\gradx u\|_{2}+ \| \dhalf u\|_{2} \qquad (u \in \mathsf{V})
\end{align}
 and implicit constants depending only upon $n$ and ellipticity constants of $A$. The reader should remark that not even the case $A^*=A$ can be treated by abstract functional analysis (while this is the case for  Kato's problem for elliptic operators) because $H$ is never self-adjoint. 
 
The original proof uses $T(b)$ technology applied with  a ``parabolic'' first order formalism of the type $DB$: $D$ is replaced by a non self-adjoint operator that involves  the non-local $\dhalf$ and it  is therefore very delicate to handle  for these two reasons. It yields the boundedness of a certain square function and, hence, the boundedness of the holomorphic functional calculus  for this first order operator. Having this at hand opens again more doors to handle boundary value problems for parabolic operators on $\R^{n+2}_{+}$.  Afterwards, a different proof of the square root estimate not using the first order approach has been recently obtained in \cite{AtEN},  even under degenerate ellipticity conditions on $A$. 

As a side product, these square functions bounds allowed us to solve the Dirichlet problem for $\partial_{t}u- \div_{\lambda,x} (A(x,t) \nabla_{\lambda,x}u)=0$ when $t\in \R, x\in \R^n, \lambda>0$ with $A$ being an $(n+1)\times (n+1)$ bounded and real matrix satisfying  \eqref{eq: Garding d} in dimension $n+1$, for data in $\L^p(\R^{n+1})$ for large $p$ \cite{AEN2}. As is in the elliptic case,  this is akin to proving that the caloric measure for this equation is $A_{\infty}$  with respect to Lebesgue measure $\d x\d t$ on $\R^{n+1}$, and this is obtained by showing a (parabolic) Carleson measure estimate for $\lambda|\nabla_{\lambda,x}u|^2\d x\d t \d \lambda$ whenever $u$ is a bounded solution.  
For this, the argument from  \cite{HKMP} in the elliptic case was streamlined to avoid some unnecessary steps and transposed to parabolic setting.   
 
\section{Conclusion}

In retrospect, although it has been abandoned later on, the  multilinear approach taken by Yves Meyer to prove the boundedness of the Cauchy integral can be considered as a miracle. Indeed,  later proofs  of that result  have been devised within complex function theory, without using multilinear analysis and real harmonic analysis.   The multilinear approach was not only useful to solve a fundamental question, but also allowed him to formulate the $T(b)$ theorem that would become a central tool. This had tremendous consequences   in singular integral operator theory and geometrical measure theory,   and also unsuspected and fruitful developments towards functional calculus and boundary value problems.  This tells us how profound is  the  impact of  Yves Meyer's original contributions on the topics evoked in this note.


\begin{thebibliography}{}

\end{thebibliography}


\begin{thebibliography}{AA}
\bibliographystyle{abbrv}

\bibitem{AAAHK} M. Alfonseca, P. Auscher, A. Axelsson, S. Hofmann,
and S. Kim, {\it Analyticity of layer potentials and $\L^{2}$ Solvability of boundary value problems for divergence form elliptic
equations with complex $L^{\infty}$ coefficients},  Adv. Math. 226 (2011), no. 5, 4533--4606.

\bibitem{AA}
A.~Amenta and P.~Auscher.
\newblock {\em Elliptic boundary value problems with fractional regularity
  data.The first order approach}, volume~37 of {\em CRM Monograph Series}.
\newblock American Mathematical Society, Providence, RI, 2018.

\bibitem{AtEN}
A. Ataei, M. Egert, K. Nystr\"om, The Kato square root problem for weighted parabolic operators, \url{arXiv:2209.11104v1}.



\bibitem{AA1}
P.~Auscher and A.~Axelsson.
\newblock Weighted maximal regularity estimates and solvability of non-smooth
  elliptic systems {I}.
\newblock {\em Invent. Math.}, 184(1) : 47--115, 2011.

\bibitem{AAH} P. Auscher, A. Axelsson, and S. Hofmann,
Functional calculus of Dirac operators and complex perturbations of
Neumann and Dirichlet problems, {\it J. Functional Analysis} {\bf 255} (2008), 374-448.

\bibitem{AAMc} P. Auscher, A. Axelsson, and A. McIntosh,
Solvability of elliptic systems with square integrable boundary data,
{\em Ark. Mat.} {\bf 48} (2010), 253--287.


\bibitem{AAMc2}
P. Auscher, A. Axelsson, and A. McIntosh,
\newblock On a quadratic estimate related to the {K}ato conjecture and boundary
  value problems.
\newblock {\em Contemp. Math. 205\/} (2010), 105--129.


\bibitem{AE-APDE}
P.~Auscher and M.~Egert.
\newblock On uniqueness results for {D}irichlet problems of elliptic systems
  without {D}e giorgi-{N}ash-{M}oser regularity.
\newblock {\em Analysis \& PDE}, 13-6 (2020), 1605--1632.

\bibitem{AE} 
P.~Auscher and M.~Egert.
\emph{Boundary value problems and Hardy spaces for elliptic systems with block structure},  	 Progress in Mathematics, 346. Birkh\"auser/Springer, Cham, 2023. xiii+310 pp.

\bibitem{AEN1} P.~Auscher, M.~Egert and K.~Nystr\"om.
$\L^2$ well-posedness of boundary value problems for parabolic systems with measurable coefficients,  {\em J. Eur. Math. Soc.}, 22(9) : 2943--3058, 2020.

\bibitem{AEN2} P.~Auscher, M.~Egert and K.~Nystr\"om.
The Dirichlet problem for parabolic operators in divergence form,   {\em J. Ecole. Polytechnique}, 5 : 407--441, 2018. 

\bibitem{AHLMcT}
P.~Auscher, S.~Hofmann, M.~Lacey, A.~McIntosh and P.~Tchamitchian.
\newblock The solution of the {K}ato square root problem for second order
  elliptic operators on {${\mathbb{R}}^n$}.
\newblock {\em Ann. of Math. (2)}, 156(2) : 633--654, 2002.

\bibitem{AHLT}  P. Auscher, 
S. Hofmann, J.L. Lewis, P. Tchamitchian,
Extrapolation of Carleson 
measures and the analyticity of 
 Kato's square root operators, 
{\it Acta Math.} {\bf 187}  (2001),  no. 2, 161--190.





\bibitem{AM}
P.~Auscher and M.~Mourgoglou.
\newblock Representation and uniqueness for boundary value elliptic problems
  via first order systems.
\newblock {\em Rev. Mat. Iberoam.}, 35(1) : 241--315, 2019.


\bibitem{AusSta}
P.~Auscher and S.~Stahlhut.
\newblock Functional calculus for first order systems of {D}irac type and
  boundary value problems.
\newblock {\em M\'{e}m. Soc. Math. Fr. (N.S.)}, (144), 2016.

\bibitem{AMcN1}
P.~Auscher, A.~M$^{\rm c}$Intosh, and A.~Nahmod.\newblock {Holomorphic
functional calculi of operators, quadratic estimates and  interpolation}.
\newblock {\em Indiana Univ. Math. J.}, 46:375--403, 1997.


\bibitem{AT1}
P.~Auscher and Ph. Tchamitchian.
\newblock {Conjecture de Kato sur les ouverts de $\R$}.
\newblock {\em Rev. Mat. Iberoamericana}, 8:149--199, 1992.

\bibitem{AT}
P.~Auscher and P.~Tchamitchian.
\newblock Calcul fonctionnel pr\'{e}cis\'{e} pour des op\'{e}rateurs elliptiques
  complexes en dimension un (et applications \`a certaines \'{e}quations
  elliptiques complexes en dimension deux).
\newblock {\em Ann. Inst. Fourier (Grenoble)}, 45(3) : 721--778, 1995.

\bibitem{AT2} P.~Auscher and P.~Tchamitchian, 
\newblock{\it
Square root problem for divergence operators and
related topics}, 
Ast\'erisque Vol. 249 (1998),
Soci\'et\'e Math\'ematique de France.

 \bibitem{ATLip}
P.~Auscher and Ph. Tchamitchian.
The square root problem for second order operators on Lipschitz domains: $\L^2$ theory, {\sl  J. Ana. Math.} {\bf 90} (2003), 1--12


\bibitem{AKMc1}
A.~Axelsson, S.~Keith and A.~McIntosh.
\newblock Quadratic estimates and functional calculi of perturbed {D}irac
  operators.
\newblock {\em Invent. Math.}, 163(3) : 455--497, 2006.

\bibitem{AKMc2}
A.~Axelsson, S.~Keith and A.~McIntosh.
The Kato square root problem for mixed boundary value problems. J. London Math. Soc. (2) 74 (2006), no. 1, 113?130.

\bibitem{CFK} L. Caffarelli, E. Fabes, C. Kenig, {\it Completely singular elliptic-harmonic measures, }Indiana Univ.
Math. J. 30 (1981), no. 6, 917--924.

\bibitem{BEH} 
S.~Bechtel, M.~Egert and R.~Haller-Dintelmann. The Kato square root problem on locally uniform domains, 
 	 {\em Adv. Math. 375} (2020), 107410, 37 pp.

\bibitem{CNS}
{A.~Castro}, {K.~Nystr\"om}, and {O.~Sande}.
\newblock {\em Boundedness of single layer potentials associated to divergence form parabolic equations with complex coefficients\/}.
\newblock Calc. Var. Partial Differ. Equ. \textbf{55} (2016), no.~5, Article:124.

\bibitem{Ch} M. Christ, A T(b) theorem with remarks on analytic 
capacity and the Cauchy integral, {\it Colloquium Mathematicum} LX/LXI 
(1990) 601-628.

\bibitem{CDM} R. Coifman, D. Deng, and Y. Meyer,
{\it 
Domaine de la racine carr\'ee de certains op\'erateurs
diff\'erentiels 
accr\'etifs},
Annales de l'Institut Fourier {\bf 33} (1983), 123--134.

\bibitem{CMcM} 
R.~Coifman, A.~McIntosh, and Y.~Meyer.
\newblock {L'int\'egrale de Cauchy d\'efinit un op\'erateur born\'e sur
  $\L^2({\R})$ pour les courbes lipschitziennes}.
\newblock {\em Ann. Math.}, 116:361--387, 1982.











\bibitem{CMast} R. 
Coifman and Y. Meyer
Au del\`a des op\'erateurs pseudo-diff\'erentiels.  Ast\'erisque, 57. Soc. Math.  France, Paris, 1978. i+185 pp.

\bibitem{CM} R. 
Coifman and Y. Meyer, ``Non-linear harmonic analysis, operator
theory, and 
PDE," in {\it Beijing Lectures in Harmonic
Analysis,} Annals of Math. Studies 
{\bf 112}, Princeton
Univ. Press, Princeton, NJ, 1986, pp. 3-45.

\bibitem{CUR}
D.~Cruz-Uribe and C.~Rios,  
\newblock{The Kato problem for operators with weighted ellipticity.}
\newblock{\em Trans. Amer. Math. Soc. 367} (2015), no. 7, 4727--4756.

\bibitem{Da}
B.E.J.~Dahlberg.
\newblock Estimates of harmonic measure.
\newblock {\em Arch. Rational Mech. Anal.}, 65(3) : 275--288, 1977.


\bibitem{Dav} G.~David.
  Analytic capacity, Calder\'on-Zygmund operators, and rectifiability. Publ. Mat. 43 (1999), no. 1, 3--25.


\bibitem{DJ}
G.~David and J.-L. Journ\'e.
\newblock {A boundedness criterion for generalized Calder\'on-Zygmund
  operators}.
\newblock {\em Ann. Math.}, 120:371--398, 1984.

\bibitem{DJS}
G.~David, J.-L. Journ\'e, and S.~Semmes.
\newblock {Op\'erateurs de Calder\'on-Zyg\-mund, fonctions para-accr\'etives et
  interpolation}.
\newblock {\em Rev. Mat. Iberoamericana}, 1:1--56, 1985.




\bibitem{FJK1} E. Fabes, D. Jerison and C. Kenig, {\it Multilinear 
square
functions and partial differential equations}, Amer. J. of Math. {\bf 
107}
(1985),
1325-1367.

 


\bibitem{H2} S. Hofmann, Local $Tb$ Theorems and applications in PDE, {\it Proceedings of the ICM Madrid}, Vol. {\bf II}, pp. 1375-1392, European Math. Soc., 2006.

\bibitem{HKMP}
{S.~Hofmann}, {C.~Kenig}, {S.~Mayboroda} and {J.~Pipher}.
\newblock {\em Square function/non-tangential maximal function estimates and the Dirichlet problem for non-symmetric elliptic operators}.
J. Amer. Math. Soc. \textbf{28} (2015), 483--529.

\bibitem{HK} S. Hofmann and S. Kim, The Green function estimates for strongly elliptic systems of second order,  {\it Manuscripta Mathematica} {\bf 124} (2007), 139-172.

\bibitem{HLMc} 
S.~Hofmann, M. Lacey and A.~M$^{\rm c}$Intosh.
\newblock {The solution of the Kato problem
for divergence form elliptic operators with Gaussian heat kernel bounds}.
\newblock {\it Annals of Math.} {\bf 156} (2002), pp 623-631.

\bibitem{HMc}  S. Hofmann and A. McIntosh, The solution of the Kato problem in two dimension,
Proceedings of the Conference on Harmonic Analysis and 
PDE held in El Escorial, Spain in July 2000, {\it Publ. Mat.} Vol. extra, 2002
pp. 143-160.

\bibitem{HN} T.~Hyt\"onen and F.~Nazarov.
The local Tb theorem with rough test functions. 
Adv. Math. 372 (2020), 107306, 36 pp. 

\bibitem{JK1} D. Jerison and C. Kenig,  The Dirichlet problem in nonsmooth domains, {\it Ann. of Math. (2)} {\bf 113}  (1981), no. 2, 367--382.

\bibitem{JK2} D. Jerison and C. Kenig,  The Neumann problem on Lipschitz domains, {\it Bull. Amer. Math. Soc. (N.S.)}  {\bf 4}  (1981), no. 2, 203--207.


\bibitem{J} J.-L. 
Journ\'e, {\it Remarks on the square root problem}, Pub.
Math. {\bf 35} 
(1991), 299-321.

\bibitem{K1} T. Kato, 
\newblock {\it Fractional powers of dissipative operators},
\newblock  J. Math. Soc.
Japan {\bf 13} (1961),  246--274. 

\bibitem{K}
T.~Kato.
\newblock {\em Perturbation theory for linear operators}.
\newblock Springer Verlag, New York, 1966.






\bibitem{KM} C. Kenig and Y. Meyer, 
{\it The Cauchy integral on Lipschitz
curves and the square root of second 
order accretive operators are the
same},
in Recent Progress in Fourier 
Analysis, I. Peral, ed., North Holland Math
Studies {\bf 111} (1985), 123-145.



 


\bibitem{Ke}
C.~E. Kenig.
\newblock {\em Harmonic analysis techniques for second order elliptic boundary
  value problems}, volume~83 of {\em CBMS Regional Conference Series in
  Mathematics}.
\newblock 
American Mathematical Society, Providence, RI, 1994.

\bibitem{KKPT}  C. Kenig, H. Koch, J. Pipher and T. Toro,  A new approach to absolute continuity of elliptic measure, with applications to non-symmetric equations. {\it Adv. Math.} {\bf 153}  (2000),  no. 2, 231--298.

\bibitem{KP} C. Kenig and J. Pipher,  The Neumann problem for elliptic equations with nonsmooth coefficients, {\it Invent. Math.} {\bf 113}  (1993),  no. 3, 447--509.

\bibitem{Li}
J.-L. Lions.
\newblock {Espaces d'interpolation et domaines de puissances fractionnaires}.
\newblock {\em J. Math. Soc. Japan}, 14:233--241, 1962.

  

\bibitem{Mc2} A. McIntosh, {\it On the Comparability of $A^{1/2}$ and
$A^{*1/2}$}, Proc. Amer. Math. Soc. 32 (1972), 430-434.


  \bibitem{Mc1} A.McIntosh, {\it Square roots of operators and applications to
hyperbolic $PDE$}, in the proceedings of the Miniconference on Operator
Theory 
and PDE, CMA, The Australian National University, Canberra, 1983.


\bibitem{Mc}
A.~McIntosh.
\newblock {Operators which have an $H^\infty$ functional calculus}.
\newblock In {\em {Miniconference on operator theory and partial differential
  equations}}, volume~14 of {\em Center for Math. and Appl.}, pages 210--231,
  Canberra, 1986. Australian National Univ.


\bibitem{Mc4} A. McIntosh, {\it The square root problem for 
elliptic
operators}, Functional Analytic Methods for
Partial Differential 
Equations, Lecture Notes in Math. vol. 1450,
Springer Verlag, 1990, pp. 
122-140.

\bibitem{McM}
A.~McIntosh and Y.~Meyer.
\newblock {Alg\`ebres d'op\'erateurs d\'efinis par des int\'e\-gra\-les
  singuli\`eres}.
\newblock {\em C. R. Acad. Sci. Paris}, 301, S\'erie 1:395--397, 1985.
  
\bibitem{Me} Y.~Meyer. 
\newblock   {Complex analysis and operator theory in Alberto Calder\'{o}n's work},
 \newblock{  
     {\em Selected papers of Alberto P. Calder\'{o}n},
   Amer. Math. Soc., Providence, RI,
  {2008},
   {593--606}.}
   
   \bibitem{N1}
\textsc{K.~Nystr\"om}.
\newblock {\em Square functions estimates and the Kato problem for second order parabolic operators in $\mathbb R^{n+1}$\/}.
\newblock Adv. Math. \textbf{293} (2016), 1--36.

\bibitem{N2}
\textsc{K.~Nystr\"om}.
\newblock {\em {$\L^2$} Solvability of boundary value problems for divergence form parabolic equations with complex coefficients\/}.
\newblock J. Differential Equations \textbf{262} (2016), no.~3, 2808--2939.

\bibitem{R}
A.~Ros\'{e}n.
\newblock Layer potentials beyond singular integral operators.
\newblock {\em Publ. Mat.}, 57(2) : 429--454, 2013.



\bibitem{S}
S.~Semmes.
\newblock {Square function estimates and the $T(b)$ Theorem}.
\newblock {\em Proc. Amer. Math. Soc.}, 110(3):721--726, 1990.


\bibitem{steinweiss} E.M.~Stein and G.~Weiss.
\newblock On the theory of $H^p$
spaces, 
\newblock {\it Acta. Math.}, 103 : 25--62, 1960.

\bibitem{T} X.~Tolsa.
Painlev\'e's problem and the semiadditivity of analytic capacity. Acta Math. 190 (2003), no. 1, 105--149.


\bibitem{V} G. Verchota, Layer potentials and regularity for the Dirichlet problem for Laplace's equation in Lipschitz domains. {\it J. Funct. Anal.} {\bf 59}  (1984),  no. 3, 572--611.

\end{thebibliography}

\end{document}